\documentclass[twoside]{article}
\usepackage{graphics}
\usepackage{graphicx}
\usepackage{latexsym}
\usepackage{amssymb}
\usepackage{amsmath}
\usepackage{BSPMsample}
\numberwithin{equation}{section}

\begin{document}

\title{Strong Forms of Weakly $e$-continuous Functions \thanks{This study was supported by Turkish-German University Scientific Research Projects Commission under the grant no: 2021BR01.}}

\author[B. S. Ayhan]{Burcu Sünbül Ayhan}

\address{Burcu Sünbül Ayhan, \\ Quality Management Unit, \\
Rectorate, Turkish-German University, \\ 34820 Beykoz-İstanbul, Turkey.}
\email{burcu.ayhan@tau.edu.tr}

\maketitle

\begin{abstract}
The main purpose of this study is to introduce and study two new classes of continuity called $eR$-continuous functions and weakly $eR$-continuous functions via $e$-regular sets. Both of the forms of continuous functions we have described are stronger than the weakly $e$-continuity. Furthermore, we obtain various characterizations of weakly $eR$-continuous functions. In addition, we examine not only the relations of these functions with some other forms of existing continuous functions, but also some of their fundamental properties.
\end{abstract}


\tableofcontents

\2010mathclass{54C08, 54D05, 54D10, 54D30.}

\section{Introduction}
There is no doubt that one of the fundamental concepts of general topology is different forms of the notion of open set. The discussion about $e$-open set types, one of the generalized open set concepts, is still a rich field to study in terms of general topology. Some forms of this concept such as $ge \Lambda$-closed sets \cite{4}, generalized $e$-closed sets \cite{1}, $\pi ge$-closed sets \cite{2}, $\Lambda_e$-sets and $V_e$-sets \cite{3} have been investigated in the recent years. Apart from these, some forms of $e$-continuity and $e$-openness of functions have been studied in \cite{a,b} as well.

In recent years, many authors have studied on generalizations of strong continuity such as strongly $\theta$-continuous functions \cite{noiri}, strongly $\theta$-precontinuous functions \cite{noiri2}, strongly $\theta$-semi continuous functions \cite{jafari-noiri}, strongly $\theta$-$b$-continuous functions \cite{park}, strongly $\theta$-$e$-continuous functions \cite{ozkoc-aslim}. On the other hand, many researchers have introduced and investigated some properties of the weakly clopen functions \cite{son-etc} and weakly $e$-continuous functions \cite{ozkoc-aslim2}.

In this paper, we investigate different classes of continuity called $eR$-continuous functions, weakly $eR$-continuous functions and strongly-$\theta$-$e$-continuous functions and study some of their fundamental properties. So, it turns out that weakly $eR$-continuous functions are weaker than strongly $\theta$-$e$-continuous functions, weakly clopen functions and $eR$-continuous functions and also stronger than weakly $e$-continuous functions.

\section{Preliminaries}
Throughout this paper, $X$ and $Y$ represent topological spaces. For a subset $A$ of a space $X$, $cl(A)$ and $int(A)$ denote the closure of $A$ and the interior of $A$, respectively. The family of every closed (resp. open, clopen) sets of $X$ is denoted by $C(X)$ (resp. $O(X), \ CO(X)$). A subset $A$ is called regular open \cite{stone} (resp. regular closed \cite{stone}) if $A=$ $int(cl(A))$ $($resp. $A=cl(int(A)))$. A point $x\in X$ is called $\delta $-cluster point \cite{velicko} of $A$ if $int(cl(U))\cap A\neq \emptyset $ for every open neighborhood $U$ of $x$. The set of all $\delta $-cluster points of $A$ is called the $\delta $-closure \cite{velicko} of $A$ and is denoted by $cl_{\delta }(A)$. If $A=cl_{\delta }(A)$, then $A$ is called $\delta $-closed \cite{velicko} and the complement of a $\delta $-closed set is called $\delta $-open \cite{velicko}. The set $\{x|(\exists U \in O(X,x))(int(cl(U)) \subseteq A)\}$ is called the $\delta$-interior of $A$ and is denoted by $int_{\delta}(A)$.

A subset $A$ of a space $X$ is called semiopen \cite{levine} (resp. preopen \cite{mashhour-monsef-deeb}, $b$-open \cite{andrijevic}, $e$-open \cite{ekici2}, $a$-open \cite{ekici1}) if $A$ $\subseteq cl(int(A))$ (resp. $A$ $\subseteq int(cl(A))$, $A$ $\subseteq cl(int(A))\cup int(cl(A))$, $A$ $\subseteq cl(int_{\delta }(A))\cup int(cl_{\delta }(A))$, $A \subseteq int(cl(int_{\delta }(A)))$. The complement of a semiopen (resp. preopen, $b$-open, $e$-open, $a$-open) set is called semiclosed \cite{levine} (resp. preclosed \cite{mashhour-monsef-deeb}, $b$-closed \cite{andrijevic}, $e$-closed \cite{ekici2}, $a$-closed \cite{ekici1}). The intersection of all semiclosed (resp. preclosed, $b$-closed, $e$-closed, $a$-closed) sets of $X$ containing $A$ is called the semi-closure \cite{levine} (resp. pre-closure \cite{mashhour-monsef-deeb}, $b$-closure \cite{andrijevic}, $e$-closure \cite{ekici2}, $a$-closure \cite{ekici1}) of $A$ and is denoted by $scl(A)$ (resp. $pcl(A)$, $bcl(A)$, $e$-$cl(A)$, $a$-$cl(A)$). The union of every semiopen (resp. preopen, $b$-open, $e$-open, $a$-open) sets of $X$ contained in $A$ is called the semi-interior \cite{levine} (resp. pre-interior \cite{mashhour-monsef-deeb}, $b$-interior \cite{andrijevic}, $e$-interior \cite{ekici2}, $a$-interior \cite{ekici1}) of $A$ and is denoted by $sint(A)$ (resp. $pint(A)$, $bint(A)$, $e$-$int(A)$, $a$-$int(A)$).

A point $x$ of $X$ is said to be $\theta $-cluster ($e$-$\theta $-cluster) point of $A$ if $cl(U)\cap A\neq \emptyset$ ($e$-$cl(U)\cap A\neq \emptyset$) for all open ($e$-open) set $U$ containing $x$. The set of all $\theta $-cluster ($e$-$\theta $-cluster) points of $A$ is called the $\theta $-closure \cite{velicko} ($e$-$\theta $-closure \cite{ozkoc-aslim}) of $A$ and is denoted by $cl_{\theta }(A)$ ($e$-$cl_{\theta }(A)$). A subset $A$ is called to be $\theta $-closed ($e$-$\theta $-closed) if $A=cl_{\theta }(A)$ ($A=e$-$cl_{\theta }(A)$). The complement of a $\theta $-closed ($e$-$\theta $-closed) set is called a $\theta$-open \cite{velicko} ($e$-$\theta $-open \cite{ozkoc-aslim}). A point $x$ of $X$ called to be a $\theta$-interior \cite{velicko} ($e$-$\theta$-interior \cite{ozkoc-aslim}) point of a subset $A$, denoted by $int_{\theta }(A)$ ($e$-$int_{\theta }(A)$), if there exists an open ($e$-open) set $U$ of $X$ containing $x$ such that $cl(U)\subseteq A$ ($e$-$cl(U)\subseteq A$).

A subset $A$ is called $e$-regular \cite{ozkoc-aslim} if it is $e$-open and $e$-closed. Also, it is noted in \cite{jumaili-yang} that 
\[ 
    e\mbox{-regular}\Rightarrow e\mbox{-}\theta \mbox{-open}\Rightarrow e\mbox{-open.}
\]
The family of every $e$-$\theta$-open $($resp. $e$-$\theta$-closed, $e$-regular, regular open, regular closed, $\delta$-open, $\delta$-closed, $\theta$-open, $\theta$-closed, semiopen, semiclosed, preopen, preclosed, $b$-open, $b$-closed, $e$-open, $e$-closed, $a$-open, $a$-closed$)$ subsets of $X$ is denoted by $e \theta O(X)$  $($resp. $e \theta C(X)$, $eR(X)$, $RO(X)$, $RC(X)$, $\delta O(X)$, $\delta C(X)$, $\theta O(X)$, $\theta C(X)$, $SO(X)$, $SC(X)$, $PO(X)$, $PC(X)$, $BO(X)$, $BC(X)$, $eO(X)$, $eC(X)$, $aO(X)$, $aC(X))$. The family of every open $($resp. closed, $e$-$\theta$-open, $e$-$\theta$-closed, $e$-regular, regular open, regular closed, $\delta$-open, $\delta$-closed, $\theta$-open, $\theta$-closed, semiopen, semiclosed, preopen, preclosed, $b$-open, $b$-closed, $e$-open, $e$-closed, $a$-open, $a$-closed$)$ sets of $X$ containing a point $x$ of $X$ is denoted by $O(X,x)$ $($resp. $C(X,x)$, $e\theta O(X,x)$, $e\theta C(X,x)$, $eR(X,x)$, $RO(X,x)$, $RC(X,x)$, $\delta O(X,x)$, $\delta C(X,x)$, $\theta O(X,x)$, $\theta C(X,x)$, $SO(X,x)$, $SC(X,x)$, $PO(X,x)$, $PC(X,x)$, $BO(X,x)$, $BC(X,x)$, $eO(X,x)$, $eC(X,x)$, $aO(X,x)$, $aC(X,x))$.

We shall use the well-known accepted language almost in the whole of the proofs of the theorems in this article. The following fundamental properties of $e$-$\theta$-closure is useful in the sequel: 

\begin{lem} \label{*} \cite{ozkoc-aslim}; \cite{jumaili-yang}
Let $A$ and $B$ be subsets of a space $X$. Then the followings are hold:\newline
\textit{(1)} $A \subseteq e$-$cl(A)$ $\subseteq$ $e$-$cl_{\theta}(A)$,\\
\textit{(2)} $e$-$cl_{\theta }(A)$ is $e$-$\theta $-closed,\newline
\textit{(3)} If $A$ is $e$-$\theta$-closed, then $A=e$-$cl_{\theta}(A)$, \newline 
\textit{(4)} If $A \subseteq B$, then $e$-$cl_{\theta}(A) \subseteq e$-$cl_{\theta}(B),$\\
\textit{(5)} $e$-$cl_{\theta}(e$-$cl_{\theta}(A))=e$-$cl_{\theta}(A)$,\\
\textit{(6)} $e$-$cl_{\theta}(X\setminus A)$ $=X\setminus e$-$int_{\theta}(A)$,\\
\textit{(7)} $x\in e$-$cl_{\theta}(A)$ iff $A\cap U\neq\emptyset$ for all $eR(X,x)$,\newline
\textit{(8)} $e$-$cl_{\theta }(A)=\bigcap \{V|(A\subseteq V)(V\in eR(X))\}=\bigcap \{V|(A\subseteq V)(V\in e \theta C(X))\}$,\newline
\textit{(9)} $A \in e \theta O(X)$ iff for all $x\in A$ there exists $U \in eR(X,x)$ such that $U\subseteq A$,\newline
\textit{(10)} Any intersection (union) of $e$-$\theta $-closed ($e$-$\theta $-open) sets is $e$-$\theta $-closed ($e$-$\theta $-open).
\end{lem}

\begin{lem} \label{bsa} \cite{ayhan-ozkoc2}
Let $A$ be a subset of a space $X$. If $A$ is an open set in $X$, then $cl(A)= cl_\theta(A)$.
\end{lem}

\begin{Defn} \label{2}
A function $f : X \rightarrow Y$ is called to be:\newline
\textit{(a)} $e$-continuous (briefly $e$.c.) \cite{ekici2} if the inverse image of all open set is $e$-open.\\
\textit{(b)} strongly $\theta$-continuous (briefly, st.$\theta$.c.) \cite{noiri} if for all $x \in X$ and all open set $V$ of $Y$ containing $f(x)$, there exists $U \in O(X,x)$ such that $f[cl(U)] \subseteq V$. \\
\textit{(c)} strongly $\theta$-semicontinuous (briefly, st.$\theta$.s.c.) \cite{jafari-noiri} if for all $x \in X$ and all open set $V$ of $Y$ containing $f(x)$, there exists $U \in SO(X,x)$ such that $f[scl(U)] \subseteq V$. \\
\textit{(d)} strongly $\theta$-precontinuous (briefly, st.$\theta$.p.c.) \cite{noiri2} if for all $x \in X$ and all open set $V$ of $Y$ containing $f(x)$, there exists $U \in PO(X,x)$ such that $f[pcl(U)] \subseteq V$. \\
\textit{(e)} strongly $\theta$-$b$-continuous (briefly, st.$\theta$.b.c.) \cite{park} if for all $x \in X$ and all open set $V$ of $Y$ containing $f(x)$, there exists $U \in BO(X,x)$ such that $f[bcl(U)] \subseteq V$. \\
\textit{(f)} strongly $\theta$-$e$-continuous (briefly, st.$\theta$.e.c.) \cite{ozkoc-aslim} if for all $x \in X$ and all open set $V$ of $Y$ containing $f(x)$, there exists $U \in eO(X,x)$ such that $ f[e\mbox{-}cl(U)] \subseteq V$.\\
\textit{(g)} weakly clopen (briefly, w.co.) \cite{son-etc} if for all $x \in X$ and all open set $V$ of $Y$ containing $f(x)$, there exists $U \in CO(X,x)$ such that $f[U] \subseteq cl(V)$.\\
\textit{(h)} weakly $b$-continuous (briefly w.$b$.c.) \cite{sengul} if for all $x \in X$ and all open set $V$ of $Y$ containing $f(x)$, there exists $U \in BO(X,x)$ such that $f[U] \subseteq cl(V)$.\\
\textit{(i)} weakly $a$-continuous (briefly w.$a$.c.) \cite{ayhan-ozkoc} if for all $x \in X$ and all open set $V$ of $Y$ containing $f(x)$, there exists $U \in aO(X,x)$ such that $f[U] \subseteq cl(V)$.\\
\textit{(j)} weakly $e$-continuous (briefly w.$e$.c.) \cite{ozkoc-aslim2} if for all $x \in X$ and all open set $V$ of $Y$ containing $f(x)$, there exists $U \in eO(X,x)$ such that $f[U] \subseteq cl(V)$.\\
\textit{(k)} weakly $BR$-continuous (briefly, w.$BR$.c.) \cite{ekicii} if for all $x \in X$ and all open set $V$ of $Y$ containing $f(x)$, there exists $U \in BR(X,x)$ such that $f[U] \subseteq cl(V)$.\\
\textit{(l)} $BR$-continuous (briefly, $BR$.c.) \cite{ekicii} if $f^{-1}[V]$ is $b$-regular in $X$ for all open set $V$ of $Y$.
\end{Defn}

%%%%%%%%%%%%%%%%%%%%%%%%%%%%%%%%%%%%%%%%%%%%%%%%%%%%%%%%%%%%%%%%%
\section{$eR$-continuity and Weakly $eR$-continuity}

\begin{Defn} \label{1}
A function $f : X \rightarrow Y$ is called weakly $eR$-continuous (briefly w.$eR$.c.) at $x\in X$ if for each open set $V$ containing $f(x)$, there exists an $e$-regular set $U$ in $X$ containing $x$ such that $f\left[ U \right] \subseteq cl(V)$. The function $f$ is w.$eR$.c. if and only if $f$ is w.$eR$.c. for all $x\in X$.
\end{Defn}

\begin{Defn} \label{11}
A function $f : X \rightarrow Y$ is called $eR$-continuous (briefly $eR$.c.) 
if $f^{-1}(V)$ is $e$-regular in $X$ for every open set $V$ of $Y$.
\end{Defn}

\begin{thm}
Let $f : X \rightarrow Y$ be a function. If $f$ is $eR$-continuous, then $f$ is weakly $eR$-continuous.
\end{thm}

\begin{pf}
Let $x \in X$ and $V \in O(Y,f(x))$.\\
$	\left.
	\begin{array}{r}
	(x \in X)(V \in O(Y,f(x))) \\
    f \text{ is $eR$.c.}
	\end{array}
	\right\}\Rightarrow  
	\begin{array}{c}
	\mbox{}\\
	\left.
	\begin{array}{r}
	\!\!\!\! f^{-1}[V] \in eR(X,x)   \\
	U:=f^{-1}[V]
	\end{array}
	\right\} \Rightarrow (U \in eR(X,x))(f[U] \subseteq V \subseteq cl(V)).
	\end{array}
$
\end{pf}

\begin{thm}
Let $f : X \rightarrow Y$ be a function. If $f$ is $eR$-continuous, then $f$ is $e$-continuous.
\end{thm}

\begin{pf}
It is obvious from Definitions \ref{2} and \ref{11}.
\end{pf}

\begin{Rem} \label{q}
From Definitions \ref{2}, \ref{1} and \ref{11}, we have the following diagram. The converses of the below implications are not true in general, as shown in the related articles and following examples.\\
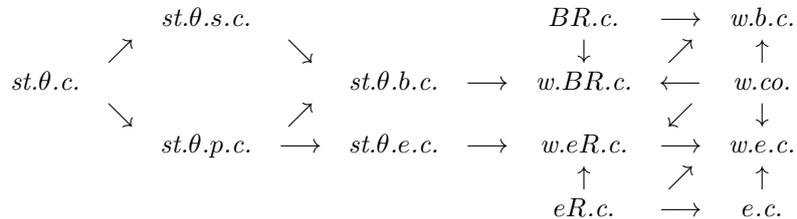
\[
\begin{array}{ccccccccc}
&  & \text{st.}\theta \text{.}s\text{.c.} &  &  &  & BR\text{.c.} & \longrightarrow  & \text{w.$b$.c.}  \\ 
& \nearrow  &  & \searrow  &  &  & \downarrow  &  \nearrow & \uparrow \\ 
\text{st.}\theta \text{.c.} &  &  &  & \text{st.}\theta \text{.}b\text{.c.} & \longrightarrow  & \text{w.}BR\text{.c.} & \longleftarrow  & \text{w.co.} \\ 
& \searrow  &  & \nearrow  &  &  &  & \swarrow  & \downarrow  \\ 
&  & \text{st.}\theta \text{.}p\text{.c.} & \longrightarrow  & \text{st.}\theta \text{.}e\text{.c.} & \longrightarrow  & \text{w.}eR\text{.c.} & \longrightarrow & \text{w.$e$.c.} \\ 
&  &  &  &  &  & \uparrow  & \nearrow  & \uparrow  \\ 
&  &  &  &  &  & eR\text{.c.} & \longrightarrow &  \text{$e$.c.} 
\end{array}%
\]
\end{Rem}

\begin{exmp}
Let $X = \{a, b, c, d, e\}$ and $\tau = \{\emptyset,\{a\},\{c\},\{a, c\},\{c, d\},\{a, c, d\},X \}$ and $\sigma = \{\emptyset, X, \{b, c, d\}\}$. Then the function $f : (X, \tau) \rightarrow (X, \sigma)$ by $f=\{(a,e),(b,b),(c,c),(d,d),(e,a)\}$ is weakly $eR$-continuous but it is not strongly $\theta$-$e$-continuous.
\end{exmp}

\begin{exmp}
Let $X = \{a, b, c, d, e\}$ and $\tau = \{\emptyset,\{a\},\{c\},\{a, c\},\{c, d\},\{a, c, d\},X \}$ and $\sigma = \{ \emptyset,X, \{b, c, d\}\}$. Then the function $f : (X, \tau) \rightarrow (X, \sigma)$ by $f=\{(a,b),(b,a),(c,c),(d,d),(e,e)\}$ is weakly $eR$-continuous but it is not $eR$-continuous.
\end{exmp}

\begin{exmp}
Let $X = \{a, b, c, d\}$ and $\tau = \{\emptyset,\{a\},\{b\},\{a, b\},\{a, b, c\},X\}$ and $\sigma = \{\emptyset, X, \{a\},\{b,c\}, \{a, b, c\}\}$. Then the function $f : (X, \tau) \rightarrow (X, \sigma)$ by $f=\{(a,a),(b,c),(c,b),(d,d)\}$ is both $e$-continuous and weakly $e$-continuous  but it is not $eR$-continuous.
\end{exmp}

\!\!\!\!\!\!\!\!\!\!\! \textbf{Question:}
Is there any weakly $e$-continuous function which is not weakly $eR$-continuous?   

\begin{thm} \label{3}
For a function $f:X\rightarrow Y$, the followings are equivalent:\newline
\textit{(1)} $f$ is weakly $eR$-continuous$;$\newline
\textit{(2)} for all $x\in X$ and all open set $V$ of $Y$ containing $f(x)$, there exists an $e$-$\theta$-open set $U$ in $X$ containing $x$ such that $f[U]\subseteq cl(V);$ \newline
\textit{(3)} $e$-$cl_{\theta }(f^{-1}[U])\subseteq f^{-1}[cl(U)]$ for all preopen set $U$ of $Y;$\newline
\textit{(4)} $f^{-1}[U]\subseteq $ $e$-$int_{\theta }(f^{-1}[cl(U)])$ for all preopen set $U$ of $Y;$\newline
\textit{(5)} $e$-$cl_{\theta }(f^{-1}[int(cl(B))])\subseteq f^{-1}[cl(B)]$ for all subset $B$ of $Y;$\newline
\textit{(6)} $e$-$cl_{\theta }(f^{-1}[int(F)])\subseteq f^{-1}[F]$ for all regular closed set $F$ of $Y;$\newline
\textit{(7)} $e$-$cl_{\theta }(f^{-1}[U])\subseteq f^{-1}[cl(U)]$ for all open subset $U$ of $Y;$\newline
\textit{(8)} $f^{-1}[U]\subseteq e$-$int_{\theta }(f^{-1}[cl(U)])$ for all open subset $U$ of $Y;$\newline
\textit{(9)} $f[e\mbox{-}cl_\theta(A)]\subseteq cl_{\theta }(f[A])$ for all subset $A$ of $X;$\newline
\textit{(10)} $e$-$cl_{\theta }(f^{-1}[B])\subseteq f^{-1}[cl_\theta(B)]$ for all subset $B$ of $Y;$\newline
\textit{(11)} $e$-$cl_{\theta }(f^{-1}[int(cl_\theta(B))])\subseteq f^{-1}[cl_\theta(B)]$ for all subset $B$ of $Y.$  
\end{thm}

\begin{pf}
    $(1) \Rightarrow (2):$ Let $x \in X$ and $V \in O(Y,f(x))$.\\
$	\left.
	\begin{array}{r}
	(x \in X)(V \in O(Y,f(x))) \\
    \mbox{ Hypothesis}
	\end{array}
	\right\}\Rightarrow  
	\begin{array}{c}
	\mbox{}\\
	\left.
	\begin{array}{r}
	\!\!\!\! (\exists U \in eR(X,x))(f[U] \subseteq cl(V))   \\
	eR(X) \subseteq e \theta O(X)
	\end{array}
	\right\} \Rightarrow 
	\end{array}
$\\
$\begin{array}{l}
\Rightarrow (\exists U \in e \theta O(X,x))(f[U] \subseteq cl(V)).
\end{array}$
\\

$(2) \Rightarrow (3):$ Let $x \in X \setminus f^{-1}[cl(G)]$ and $G \in PO(Y)$.\\
$
\left.\begin{array}{rr}
x \in X \setminus f^{-1}[cl(G)] \Rightarrow f(x) \in Y \setminus cl(G) \Rightarrow (\exists V \in  O(Y,f(x)))(V \cap G = \emptyset)\\
\mbox{Hypothesis} 
\end{array}\right\} \Rightarrow 
$\\
$\begin{array}{l}
\Rightarrow (\exists U \in e \theta O(X,x))(f[U] \subseteq cl(V))(V \cap G = \emptyset)
\end{array}$\\
$\begin{array}{l}
\Rightarrow (\exists U \in e \theta O(X,x))(f[U] \cap  G = \emptyset)
\end{array}$\\
$\begin{array}{l}
\Rightarrow (\exists U \in e \theta O(X,x))(U \cap  f^{-1}[G] = \emptyset)
\end{array}$\\
$\begin{array}{l}
\Rightarrow x \notin e\mbox{-}cl_\theta(f^{-1}[G])
\end{array}$\\
$\begin{array}{l}
\Rightarrow x \in X \setminus e\mbox{-}cl_\theta(f^{-1}[G]).
\end{array}$\\ 

$(3) \Rightarrow (4):$ Let $G \in PO(Y)$.\\
$
\left.\begin{array}{rr}
G \in PO(Y) \Rightarrow Y \setminus cl(G) \in O(Y) \subseteq PO(Y)\\
\mbox{Hypothesis} 
\end{array}\right\} \Rightarrow  e\mbox{-}cl_\theta ( f^{-1}[Y\setminus cl(G)]) \subseteq f^{-1}[cl(Y\setminus cl(G))]
$\\
$\begin{array}{l}
\Rightarrow X \setminus e\mbox{-}int_\theta ( f^{-1}[cl(G)]) \subseteq X \setminus f^{-1}[int(cl(G))] \subseteq X \setminus f^{-1}[G]
\end{array}$ \\
$\begin{array}{l}
\Rightarrow  f^{-1}[G] \subseteq  e\mbox{-}int_\theta ( f^{-1}[cl(G)]).
\end{array}$ \\ 

$(4) \Rightarrow (5):$ Let $B \subseteq Y$.\\
$
\left.\begin{array}{rr}
B \subseteq Y \Rightarrow Y \setminus cl(B) \in O(Y) \subseteq PO(Y)\\
\mbox{Hypothesis} 
\end{array}\right\} \Rightarrow f^{-1}[Y\setminus cl(B)] \subseteq e\mbox{-}int_\theta ( f^{-1}[cl(Y\setminus cl(B))])
$\\
$\begin{array}{l}
\Rightarrow X \setminus f^{-1}[cl(B)] \subseteq  X \setminus e\mbox{-}cl_\theta ( f^{-1}[int(cl(B))])
\end{array}$ \\
$\begin{array}{l}
\Rightarrow e\mbox{-}cl_\theta ( f^{-1}[int(cl(B))]) \subseteq f^{-1}[cl(B)].
\end{array}$ \\ 

$(5) \Rightarrow (6):$ Let $F \in RC(Y)$.\\
$
\left.\begin{array}{rr}
F \in RC(Y) \Rightarrow int(F) \subseteq Y\\
\mbox{Hypothesis} 
\end{array}\right\} \Rightarrow $\\
$\begin{array}{l}
\Rightarrow e\mbox{-}cl_\theta (f^{-1}[int(F)]) =  e\mbox{-}cl_\theta (f^{-1}[int(cl(int(F)))]) \subseteq f^{-1}[cl(int(F))] = f^{-1}[F].
\end{array}
$\\ 

$(6) \Rightarrow (7):$ Let $U \in O(Y)$.\\
$
\left.\begin{array}{rr}
U \in O(Y) \Rightarrow cl(U) \in RC(Y)\\
\mbox{Hypothesis} 
\end{array}\right\} \Rightarrow e\mbox{-}cl_\theta (f^{-1}[U]) \subseteq  e\mbox{-}cl_\theta (f^{-1}[int(cl(U))]) \subseteq f^{-1}[cl(U)].
$\\ 

$(7) \Rightarrow (8):$ Let $U \in O(Y)$.\\
$
\left.\begin{array}{rr}
U \in O(Y) \Rightarrow Y \setminus cl(U) \in O(Y) \\
\mbox{Hypothesis} 
\end{array}\right\} \Rightarrow e\mbox{-}cl_\theta (f^{-1}[Y \setminus cl(U)]) \subseteq f^{-1}[cl(Y \setminus cl(U))]
$\\
$\begin{array}{l}
\Rightarrow  X \setminus e\mbox{-}int_\theta (f^{-1}[cl(U)]) \subseteq X \setminus f^{-1}[int(cl(U))] \subseteq X \setminus  f^{-1}[U] 
\end{array}$ \\
$\begin{array}{l}
\Rightarrow f^{-1}[U] \subseteq f^{-1}[int(cl(U))] \subseteq e\mbox{-}int_\theta (f^{-1}[cl(U)]).
\end{array}$ \\

$(8) \Rightarrow (9):$ Let $A \subseteq X$.\\
$
\left.\begin{array}{rr}
A \subseteq X \Rightarrow int(Y \setminus f[A]) \in O(Y) \\
\mbox{Hypothesis} 
\end{array}\right\} \Rightarrow  f^{-1}[int(Y \setminus f[A])] \subseteq e\mbox{-}int_\theta (f^{-1}[cl(int(Y \setminus f[A]))])
$\\
$\begin{array}{l}
\Rightarrow  X \setminus f^{-1}[cl(f[A])] \subseteq X \setminus e\mbox{-}cl_\theta (f^{-1}[int(cl(f[A]))])
\end{array}$ \\
$\begin{array}{l}
\Rightarrow  f[e\mbox{-}cl_\theta (f^{-1}[int(cl(f[A]))])]\subseteq f[f^{-1}[cl(f[A])]] 
\end{array}$ \\
$\begin{array}{l}
\Rightarrow  f[e\mbox{-}cl_\theta(A)]\subseteq cl_{\theta }(f[A]).
\end{array}$ \\

$(9) \Rightarrow (10):$ Let $B \subseteq Y$.\\ 
$
\left.\begin{array}{rr}
B \subseteq Y \Rightarrow f^{-1}[B] \subseteq X \\
\mbox{Hypothesis} 
\end{array}\right\} \Rightarrow f[e\mbox{-}cl_\theta (f^{-1}[B])] \subseteq cl_\theta (f[f^{-1}[B]]) \subseteq cl_\theta (B)
$\\
$\begin{array}{l}
\Rightarrow e\mbox{-}cl_\theta (f^{-1}[B]) \subseteq f^{-1}[f[e\mbox{-}cl_\theta (f^{-1}[B])]] \subseteq f^{-1}[cl_\theta (f[f^{-1}[B]])] \subseteq f^{-1}[cl_\theta (B)].
\end{array}$ \\ 

$(10) \Rightarrow (11):$ Let $B \subseteq Y$.\\
$
\left.\begin{array}{rr}
B \subseteq Y \Rightarrow int(cl_\theta (B)) \subseteq Y \\
\mbox{Hypothesis} 
\end{array}\right\} \Rightarrow e\mbox{-}cl_\theta (f^{-1}[int(cl_\theta (B))]) \subseteq f^{-1}[cl_\theta (int(cl_\theta (B)))] \subseteq f^{-1}[cl_\theta (B)].
$\\ 

$(11) \Rightarrow (1):$ Let $x \in X$ and $U \in O(Y,f(x))$.\\
$	\left.
	\begin{array}{r}
	(x \in X)(U \in O(Y,f(x))) \\
    \mbox{ Hypothesis}
	\end{array}
	\right\}\Rightarrow  
	\begin{array}{c}
	\mbox{}\\
	\left.
	\begin{array}{r}
	\!\!\!\!\! e\mbox{-}cl_\theta (f^{-1}[int(cl_\theta (U))]) \subseteq  f^{-1}[cl_\theta (U)]\\
	U \subseteq int(cl(U))=int(cl_\theta(U))
	\end{array}
	\right\} \Rightarrow 
	\end{array}
$\\
$\begin{array}{l}
\Rightarrow e\mbox{-}cl_\theta (f^{-1}[U]) \subseteq e\mbox{-}cl_\theta (f^{-1}[int(cl_\theta (U))]) \subseteq f^{-1}[cl_\theta (U)]  \overset{\text{Lemma \ref{bsa}} }{=} f^{-1}[cl(U)]
\end{array}$ \\

It follows by $(1) \Leftrightarrow (7)$ that $f$ is weakly $eR$ continuous.
\end{pf}

\begin{Defn} \cite{jumaili-yang}
A function $f : X \rightarrow Y$ is called to be contra $e\theta$-continuous (briefly c.$e\theta$.c.)  if $f^{-1}[V]$ is $e$-$\theta$-closed in $X$ for every open set $V$ of $Y$.
\end{Defn}

\begin{thm}
    If $f : X \rightarrow Y$ is contra $e\theta$-continuous, then $f$ is weakly $eR$-continuous.
\end{thm}

\begin{pf}
    Let $V \in O(Y)$.\\
$
\left.\begin{array}{rr}
V \in O(Y)\\
f \mbox{ is c.$e\theta$.c.} 
\end{array}\right\} \Rightarrow f^{-1}[V] \in e\theta C(X) \Rightarrow e\mbox{-}cl_\theta(f^{-1}[V])=f^{-1}[V] \subseteq f^{-1}[cl(V)]
$\\
Then by \ref{3}(7) $f$ is weakly $eR$-continuous.
\end{pf}

\begin{thm} \label{4}
For a function $f : X \rightarrow Y$, the followings are equivalent:\newline
\textit{(1)} $f$ is weakly $eR$-continuous at $x \in X$;\newline
\textit{(2)} $x \in e\mbox{-}int_{\theta}(f^{-1}[cl(V)])$ for each neighborhood $V$ of $f(x)$.
\end{thm}

\begin{pf}
$(1) \Rightarrow (2):$ Let $V \in O(Y,f(x))$.\\
$
\left.\begin{array}{rr}
V \in O(Y,f(x))\\
\mbox{Hypothesis} 
\end{array}\right\} \Rightarrow (\exists U \in eR(X,x))(f[U] \subseteq cl(V)) \Rightarrow (\exists U \in eR(X,x))(U \subseteq f^{-1}[cl(V)])
$\\
$\begin{array}{l}
\Rightarrow (\exists U \in eR(X,x))(U \subseteq e\mbox{-}int_{\theta}(U) \subseteq e\mbox{-}int_\theta(f^{-1}[cl(V)]))
\end{array}$\\
$\begin{array}{l}
\Rightarrow x \in e\mbox{-}int_{\theta}(f^{-1}[cl(V)]).
\end{array}$\\ 

$(2) \Rightarrow (1):$ Let $V \in O(Y,f(x))$.\\
$
\left.\begin{array}{rr}
V \in O(Y,f(x))\\
\mbox{Hypothesis} 
\end{array}\right\} \Rightarrow x \in e\mbox{-}int_{\theta}(f^{-1}[cl(V)]) \Rightarrow (\exists U \in eR(X,x))(f[U] \subseteq cl(V))
$\\

Thus, $f$ is weakly $eR$-continuous at $x \in X$.
\end{pf}

\begin{lem} \cite{ozkoc-izci} \label{5}
Let $X$ be a topological space and $A,B \subseteq X$. If $A \in aO(X)$ and $B \in eO(X)$, then $A \cap B \in eO(X)$.
\end{lem}

\begin{thm}
If $f : X \rightarrow Y$ is weakly $eR$-continuous at $x \in X$, then there exists a nonempty $e$-open set $U \subseteq H$ such that $U \subseteq e\mbox{-}cl_{\theta }(f^{-1}[cl(V)])$ for every neighborhood $V$ of $f(x)$ and every $a$-open neighborhood $H$ of $x$.
\end{thm}

\begin{pf}
Let $H \in aO(X,x)$ ve $V \in O(Y,f(x))$.\\
$
\left.\begin{array}{rr}
(H \in aO(X,x))(V \in O(Y,f(x)))\\
f \mbox{ is w.$eR$.c. at $x \in X$}
\end{array}\right\} \overset{\text{Theorem } \ref{4}}{\Rightarrow } x \in e\mbox{-}int_{\theta}(f^{-1}[cl(V)])
$\\
$
\left.\begin{array}{rr}
 \overset{\text{Lemma } \ref{*}(9)}{\Rightarrow} (\exists G \in eR(X,x))(G \subseteq e\mbox{-}int_{\theta}(f^{-1}[cl(V)]))\\
U:=G \cap H 
\end{array}\right\} \overset{\text{Lemma } \ref{5}}{\Rightarrow } 
$\\
$\begin{array}{l}
\Rightarrow   (U \in eO(X,x))(U \subseteq H)(U \subseteq G \subseteq e\mbox{-}int_{\theta }(f^{-1}[cl(V)]) \subseteq e\mbox{-}cl_{\theta }(f^{-1}[cl(V)])).
\end{array}$   
\end{pf}

\begin{thm}
Let $f : X \rightarrow Y$ be a function. If $f^{-1}[cl_\theta(U)]$ is $e$-$\theta$-closed in $X$ for every subset $U$ of $Y$, then $f$ is weakly $eR$-continuous.
\end{thm}

\begin{pf}
Let $f^{-1}[cl_\theta(U)]$ is $e$-$\theta$-closed in $X$ for every subset $U$ of $Y$.\\
$\begin{array}{l}
(U \subseteq Y)(f^{-1}[cl_\theta (U)]) \in e\theta C(X)) \Rightarrow e\text{-}cl_\theta (f^{-1}[U]) \subseteq e\text{-}cl_\theta (f^{-1}[cl_\theta (U)])=f^{-1}[cl_\theta (U)]
\end{array}$  \\
Then by Theorem \ref{3}(10), $f$ is weakly $eR$-continuous.
\end{pf}

\begin{thm} \label{b}
Let $f : X \rightarrow Y$ be a weakly $eR$-continuous function. Then the followings hold:\\
\textit{(1)} $f^{-1}[U]$ is $e$-$\theta$-open in $X$ for each $\theta$-open set $U$ of $Y$,\\
\textit{(2)} $f^{-1}[V]$ is $e$-$\theta$-closed in $X$ for each $\theta$-closed set $V$ of $Y$.
\end{thm}

\begin{pf}
It follows from Theorem \ref{3}.
\end{pf}

\begin{thm} \cite{ozkoc-aslim} \label{6}
A function $f: X \rightarrow Y$ strongly $\theta$-$e$-continuous if and only if for all $x \in X$ and all open set $V$ of $Y$ containing $f(x)$, there exists $U \in eR(X,x)$ such that $f(U) \subseteq V$.
\end{thm}

\begin{thm}
Let $f: X \rightarrow Y$ be a function. If $Y$ is regular, then the followings are equivalent:\\
\textit{(1)} $f$ is weakly $eR$-continuous;\\
\textit{(2)} $f$ is strongly $\theta$-$e$-continuous.
\end{thm}

\begin{pf}
$(1) \Rightarrow (2):$ Let $x \in X$ and $V \in O(Y,f(x))$.\\
$	\left.
	\begin{array}{r}
	(x \in X)(V \in O(Y,f(x))) \\
    Y \mbox{ is regular}
	\end{array}
	\right\}\Rightarrow  
	\begin{array}{c}
	\mbox{}\\
	\left.
	\begin{array}{r}
	\!\!\!\! (\exists H \in O(Y,f(x)))(cl(H) \subseteq V)   \\
	f \mbox{ is w.$eR$.c.}
	\end{array}
	\right\} \Rightarrow 
	\end{array}$\\
$\begin{array}{l}
\Rightarrow (\exists U \in eR(X,x))(f[U] \subseteq cl(H) \subseteq V)
\end{array}$\\ 
Then by Theorem \ref{6} $f$ is strongly $\theta$-$e$-continuous.\\ 

$(2) \Rightarrow (1):$ It is obvious from Definitions \ref{2} and \ref{1}.
\end{pf}

\section{Some Basic Properties of Weakly $eR$-continuity}

\begin{Defn} 
A space $X$ is called to be $e$-connected \cite{ekici3} if $X$ is not the union of two disjoint nonempty $e$-open sets.
\end{Defn}

\begin{thm}
If $f : X \rightarrow Y$ is a weakly $eR$-continuous surjection and $X$ is $e$-connected, then $Y$ is connected.
\end{thm}

\begin{pf}
    Suppose that $Y$ is not connected.\\
$\left. 
	\begin{array}{c}
	Y\mbox{ is not connected}\Rightarrow (\exists A,B\in O(Y) \setminus
	\{\emptyset\} )(A\cap B=\emptyset )(A\cup B=Y)
	\end{array}
	\right. $\newline
	$\left. 
	\begin{array}{r}
	\Rightarrow (A,B\in CO(Y)\setminus \{\emptyset \})(A\cap B=\emptyset )(A\cup
	B=Y) \\ 
	f\mbox{ is w.$eR$.c.}
	\end{array}
	\right\} \overset{\text{Theorem \ref{3}(8)}}{\Rightarrow } $\newline
	$
	\begin{array}{r}
	\Rightarrow f^{-1}[A]\subseteq e\mbox{-}int_\theta(f^{-1}[cl(A)])=e\mbox{-}int_\theta(f^{-1}[A])
	\end{array}$\\
	$\left. 
	\begin{array}{r}
	\! \Rightarrow (f^{-1}[B]\subseteq e\mbox{-}int_\theta(f^{-1}[cl(B)])=e\mbox{-}int_\theta(f^{-1}[B]))(f^{-1}[A\cap	B]=\emptyset )(f^{-1}[A\cup B]=X) \\ 
	f\mbox{ is surjection}
	\end{array}
	\right\} \Rightarrow $\newline
	$\left. 
	\begin{array}{c}
	\Rightarrow (f^{-1}[A],f^{-1}[B]\in eO(X)\setminus \{\emptyset
	\})(f^{-1}[A]\cap f^{-1}[B] =\emptyset )(f^{-1}[A]\cup
	f^{-1}[B] =X)
	\end{array}
	\right. $ \\
 
 This is a contradiction to the fact that $X$ is $e$-connected. Thus, $Y$ is connected.
\end{pf}

\begin{Defn}
A space $X$ is called:\\
$\textit{(1)}$ Urysohn \cite{willard} if for all distinct two points $x$ and $y$ in $X$, there exist open sets $U$ and $V$ such that $x \in U$, $y \in V$ and $cl(U) \cap cl(V) = \emptyset$.\\
$\textit{(2)}$ Clopen $T_2$ \cite{ekici} if for all distinct two points $x$ and $y$ in $X$, there exist disjoint clopen sets $U$ and $V$ of $X$ such that $x \in U$, $y \in V$ and $U \cap V = \emptyset$.\\
$\textit{(3)}$ $eR$-$T_1$ if for all distinct two points $x$ and $y$ in $X$, there exist $e$-regular sets $U$ and $V$ of $X$ containing $x$ and $y$, respectively, such that $y \notin U$ and $x \notin V$.\\
$\textit{(4)}$ $eR$-$T_2$ if for all distinct two points  $x$ and $y$ in $X$, there exist disjoint $e$-regular sets $U$ and $V$ of $X$ such that $x \in U$, $y \in V$ and $U \cap V = \emptyset$.
\end{Defn}

\begin{Rem}
Every clopen $T_2$ space is $eR$-$T_2$. This implication is not reversible as shown in the following example.
\end{Rem}

\begin{exmp}
Let $X=\{a, b, c, d\}$ and $\tau=\{\emptyset,\{a\},\{b\},\{a, b\},\{a, b, c\},\{a, b, d\}, X\}$. It is not difficult to see $CO(X)=\{ \emptyset, X\}$ and $eR(X)=2^X \setminus \{\{c\},\{d\},\{c,d\}, \{a,b\}, \{a,b,c\},\{a,b,d\}\}$. Then $(X, \tau)$ is $eR$-$T_2$ but it is not clopen $T_2$.
\end{exmp}

\begin{thm}
Let $f : X \rightarrow Y$ is a weakly $eR$-continuous function and $g : X \rightarrow Y$ is a weakly $a$-continuous function. If $Y$ is Urysohn, then $A=\{x|f(x)=g(x)\} \in eC(X)$.
\end{thm}

\begin{pf}
 Let $x \notin A $.\\
        $\left. 
	\begin{array}{r}
	x \notin A  \Rightarrow f(x) \neq g(x)\\ 
	Y \text{ is Urysohn}
	\end{array}
	\right\} \Rightarrow $\\
$\left. 
	\begin{array}{r}
	\Rightarrow (\exists V_1 \in O(Y,f(x)))(\exists V_2 \in O(Y,g(x)))(cl(V_1) \cap cl(V_2)=\emptyset) \\ 
	(f \text { is w.$eR$.c.})(g \text { is w.$a$.c.})
	\end{array}
	\right\} \Rightarrow  $\\
$
\left.\begin{array}{rr}
 \Rightarrow (\exists U \in eR(X,x))(\exists G \in aO(X,x))(f[U] \subseteq cl(V_1))(g[G] \subseteq cl(V_2))\\
W:=U \cap G 
\end{array}\right\} \overset{\text{Lemma } \ref{5}}{\Rightarrow} $\\
$\begin{array}{l}
\Rightarrow (\exists W \in eO(X,x))(f[W] \cap g[W] \subseteq f[U] \cap g[G] \subseteq cl(V_1) \cap cl(V_2)=\emptyset)
\end{array}$\\
$\begin{array}{l}
\Rightarrow (\exists W \in eO(X,x))(W \cap A=\emptyset)
\end{array}$\\
$\begin{array}{l}
\Rightarrow x \notin e\text{-}cl(A).
\end{array}$
\end{pf}

\begin{thm}
Let $f: X \rightarrow Y$ be a weakly $eR$-continuous injection. If $Y$ is Hausdorff, then $X$ is $eR$-$T_1$.
\end{thm}

\begin{pf}
    Let $x,y \in X$ and $x \neq y$.\\
$\left. 
\begin{array}{r}
(x,y \in X)(x \neq y) \overset{f\text{ is injective}}{\Rightarrow} f(x) \neq f(y)\\ 
Y \mbox{ is Hausdorff}
\end{array}
\right\} \Rightarrow  (\exists U \in O(Y,f(x))(\exists V \in O(Y,f(y))(U \cap V=\emptyset)$\\
$\left. 
\begin{array}{r}
\Rightarrow (\exists U \in O(Y,f(x))(\exists V \in O(Y,f(y))(f(x) \notin cl(V))(f(y) \notin cl(U))\\ 
f \mbox{ is w.$eR$.c.}
\end{array}
\right\} \Rightarrow $\\	
$\begin{array}{l}
\Rightarrow (\exists A \in eR(X,x))(\exists B \in eR(X,y))(f[A] \subseteq  cl(U))(f[B] \subseteq  cl(V))
\end{array}$\\
$\begin{array}{l}
\Rightarrow (\exists A \in eR(X,x))(\exists B \in eR(X,y))(y \notin A)(x \notin B).
\end{array}$
\end{pf}

\begin{thm}
    Let $f: X \rightarrow Y$ be a weakly $eR$-continuous injection. If $Y$ is Urysohn, then $X$ is $eR$-$T_2$.
\end{thm}

\begin{pf}
    Let $x,y \in X$ and $x \neq y$.\\
$\left. 
\begin{array}{r}
(x,y \in X)(x \neq y) \overset{f\text{ is injective}}{\Rightarrow} f(x) \neq f(y)\\ 
Y \mbox{ is Urysohn}
\end{array}
\right\} \Rightarrow $\\
$\left. 
\begin{array}{r}
\Rightarrow (\exists U \in O(Y,f(x))(\exists V \in O(Y,f(y))(cl(U) \cap cl(V)=\emptyset) \\ 
f \mbox { is w.$eR$.c.}
\end{array}
\right\} \Rightarrow $\\
$\begin{array}{l}
\Rightarrow (\exists A \in eR(X,x))(\exists B \in eR(X,y))(f[A] \subseteq  cl(U))(f[B] \subseteq  cl(V))
\end{array}$\\
$\begin{array}{l}
\Rightarrow (\exists A \in eR(X,x))(\exists B \in eR(X,y))(A \subseteq  f^{-1}[cl(U)])(B \subseteq  f^{-1}[cl(V)])
\end{array}$\\
$\begin{array}{l}
\Rightarrow (\exists A \in eR(X,x))(\exists B \in eR(X,y))(A \cap B \subseteq  f^{-1}[cl(U)] \cap f^{-1}[cl(V)]=\emptyset).
\end{array}$
\end{pf}

\begin{cor}
Let $f: X \rightarrow Y$ be a weakly clopen injection. If $Y$ is Urysohn, then $X$ is $eR$- $T_2$. 
\end{cor}

%\begin{lem}
%Let $A$ and $B$ be subsets of a space $X$ and $Y$, respectively. Then the followings are hold:\\
%\textit{(1)} \cite{ozkoc-izci} If $A \in eO(X)$ and $B \in eO(Y)$, then $A \times B \in eO(X \times Y)$.\\
%\textit{(2)} $e\mbox{-}cl(A \times B) \subseteq e\mbox{-}cl(A) \times e\mbox{-}cl(B)$.
%\end{lem}

%\begin{pf}
%\textit{(2)} Let $(x,y) \notin e\mbox{-}cl(A) \times e\mbox{-}cl(B)$.\\
%$\begin{array}{l}
%(x,y) \notin e\mbox{-}cl(A) \times e\mbox{-}cl(B) \Rightarrow x \notin e\mbox{-}cl(A) \vee y \notin e\mbox{-}cl(B)\\
%\end{array}$\\
%$\begin{array}{l}
%\Rightarrow (\exists U_1 \in eO(X,x))(U_1 \cap A = \emptyset) \vee (\exists U_2 \in eO(Y,y))(U_2 \cap B = \emptyset)
%\end{array}$\\
%$\begin{array}{l}
%\Rightarrow (\exists U_1 \in eO(X,x))(\exists U_2 \in eO(Y,y))(U_1 \cap A = \emptyset \vee U_2 \cap B = \emptyset)
%\end{array}$\\
%$\begin{array}{l}
%\Rightarrow (\exists U_1 \in eO(X,x))(\exists U_2 \in eO(Y,y))((U_1 \cap A) \times (U_2 \cap B) = \emptyset)
%\end{array}$\\
%$\begin{array}{l}
%\Rightarrow (U_1 \times U_2 \in eO(X \times Y,(x,y)))((U_1 \times U_2) \cap (A \times B) = \emptyset)
%\end{array}$\\
%$\begin{array}{l}
%\Rightarrow (x,y) \notin e\mbox{-}cl(A \times B).
%\end{array}$
%\end{pf}

Let $\{X_{\alpha }|\alpha \in I\}$ and $\{Y_{\alpha }|\alpha \in I \}$
be any two families of topological spaces with the same index set $I.$ The
product space of $\{X_{\alpha }|\alpha \in I\}$ (resp. $\{Y_{\alpha }|\alpha
\in I\}$) is simply denoted by $\Pi X_{\alpha }$ (resp. $\Pi Y_{\alpha }$).
Let $f_{\alpha }:X_{\alpha }\rightarrow Y_{\alpha }$ be a function for all $
\alpha \in I.$ Let $f:\Pi X_{\alpha }\rightarrow \Pi Y_{\alpha }$ be the
product function defined as follows: $f(\{x_{\alpha }\})=\{f_{\alpha
}(x_{\alpha })\}$ for all $\{x_{\alpha }\}\in \Pi X_{\alpha }.$ 

\begin{thm}
    Let $\{X_{\alpha }|\alpha \in I\}$ and $\{Y_{\alpha }|\alpha \in I \}$ be any two families of topological spaces. If $f_\alpha: X_\alpha \rightarrow Y_\alpha$ is weakly $eR$-continuous for each $\alpha \in I$, then the function  $f:\Pi X_{\alpha }\rightarrow \Pi Y_{\alpha }$ is weakly $eR$-continuous.
\end{thm}

\begin{pf}
    Let $x=\{x_{\alpha }\}\in \Pi X_{\alpha }$ and $V\in O(\Pi Y_{\alpha },f(x))$.\\
        $\begin{array}{l}
        V\in O(\Pi Y_{\alpha },f(x))\Rightarrow (\exists J=\{\alpha _{1},\alpha	_{2},\ldots ,\alpha _{n}\}\subseteq I)
        \end{array}$\\
        $\left. 
	\begin{array}{r}
	\left( W_{\alpha }:=\left\{ 
	\begin{array}{ccc}
	W_{\alpha _{j}}\in O(Y_{\alpha _{j}}) & , & \alpha \in J \\ 
	Y_{\alpha } & , & \alpha \notin J%
	\end{array}
	\right. \right)\left( 
	\Pi W_{\alpha }\in O(\Pi Y_{\alpha },f(x))\right) \left( \newline
	\Pi W_{\alpha }\subseteq V\right)  \\ 
	(\forall \alpha \in I)(f_{\alpha }\mbox{ is w.$eR$.c.})
	\end{array}
	\right\} \Rightarrow $\newline
	$\left. 
	\begin{array}{r}
	\Rightarrow (\exists U_{\alpha }\in eR(X_{\alpha },x_{\alpha
	}))\left( f_{\alpha }[U_{\alpha }]\subseteq cl(W_{\alpha })\right)  \\ 
	U:=\prod\limits_{j=1}^{n}U_{\alpha _{j}}\times \prod\limits_{\alpha \notin
		J}X_{\alpha }
	\end{array}
	\right\}  \Rightarrow  $\\
	$\left. 
	\begin{array}{c}
	\Rightarrow (U\in eR(\Pi X_{\alpha },x))\left( f[U]\subseteq
	\prod\limits_{j=1}^{n}f_{\alpha }[ U_{\alpha _{j}}] \times
	\prod\limits_{\alpha \notin J}Y_{\alpha }\subseteq
	\prod\limits_{j=1}^{n}cl( W_{\alpha _{j}}) \times
	\prod\limits_{\alpha \notin J}Y_{\alpha }\subseteq cl(V)\right).
	\end{array}
	\right. $
\end{pf}

Recall that the graph of a function $f: X \rightarrow Y$ is the subset $\{(x,f(x))|x \in X\}$ of the product space $X \times Y$ and denoted by $G(f)$.

\begin{thm}
Let $f : X \rightarrow Y$ be a function. If the graph function $g$ is weakly $eR$-continuous, then $f$ is weakly $eR$-continuous.
\end{thm}

\begin{pf}
 Let $x \in X$ and $V \in O(Y,f(x))$.\\
$
\left.\begin{array}{rr}
(x \in X)(V \in O(Y,f(x))) \Rightarrow X \times V \in O(X \times Y,g(x))\\
g \mbox{ is w.$eR$.c.} 
\end{array}\right\} \Rightarrow 
$\\
$\begin{array}{l}
 \Rightarrow  (\exists U \in eR(X,x))(g[U] \subseteq cl(X \times V)= X \times cl(V))
\end{array}$  
\\
$\begin{array}{l}
 \Rightarrow  (\exists U \in eR(X,x))(f[U] \subseteq cl(V)).
\end{array}$     
\end{pf}

\begin{Defn} \label{9}
A function $f: X \rightarrow Y$ has an $er$-graph if for all $(x,y) \notin G(f)$, there exist $U \in eR(X,x)$ and $V \in O(Y,y)$ such that $(U \times cl(V)) \cap G(f)= \emptyset$.
\end{Defn}

\begin{lem} \label{10}
A function $f: X \rightarrow Y$ has an $er$-graph if and only if for all $(x,y)\notin G(f)$, there exist $U\in eR(X,x)$ and $V \in O(Y,y)$  such that $f\left[ U\right] \cap cl(V)=\emptyset .$
\end{lem}

\begin{pf}
It is obvious from Definition \ref{9}.
\end{pf}

\begin{thm}
If $f: X \rightarrow Y$ is weakly $eR$-continuous and $Y$ is a Urysohn space, then $G(f)$ is an $er$-graph.
\end{thm}

\begin{pf}
Let $(x,y) \notin G(f).$\\
$
\left.\begin{array}{rr}
(x,y) \notin G(f) \Rightarrow y \neq f(x)\\
Y \mbox{ is Urysohn}
\end{array}\right\} \Rightarrow 
$\\
$
\left.\begin{array}{rr}
\Rightarrow (\exists A\in O(X,f(x)))(\exists B\in O(Y,y))(cl(A) \cap cl(B)=\emptyset )\\
f\mbox{ is w.$eR$.c.}
\end{array}\right\} \Rightarrow 
$\\
 $\begin{array}{c}
\Rightarrow (\exists G\in eR(X,x))(\exists B\in O(Y,y))(f[G]\subseteq cl(A))
\end{array} $\\
$\begin{array}{c}
\Rightarrow (\exists G\in eR(X,x))(\exists B\in O(Y,y))(f[G] \cap cl(B)= \emptyset)
\end{array}$\\
$\begin{array}{c}
\Rightarrow (\exists G\in eR(X,x))(\exists B\in O(Y,y))((G\times cl(B))\cap G(f)=\emptyset).
\end{array}$
\end{pf}

\begin{thm}
If $f: X \rightarrow Y$ has an $er$-graph and a weakly $eR$-continuous injection, then $X$ is $eR$-$T_2$.  
\end{thm}

\begin{pf}
Let $x,y \in X$ and $x \neq y$.\\
$
\left.\begin{array}{rr}
(x,y \in X)(x \neq y)  \overset{f \text{ is injective}}{\Rightarrow} f(x) \neq f(y)  \Rightarrow  (x,f(y)) \notin G(f) \\
G(f) \mbox{ is $er$-graph}
\end{array}\right\} \Rightarrow $\\
$\begin{array}{l}
\Rightarrow (\exists U\in eR(X,x))(\exists V\in O(Y,f(y)))( (U\times cl(V))\cap G(f)=\emptyset) 
\end{array}$\\
$
\left.\begin{array}{rr}
\Rightarrow (\exists U\in eR(X,x))(\exists V\in O(Y,f(y)))( f[U]\cap cl(V)=\emptyset)   \\
f \mbox{ is w.$eR$.c.}
\end{array}\right\} \Rightarrow $\\
 $\begin{array}{c}
\Rightarrow (\exists U\in eR(X,x))(\exists W\in eR(X,y))(f[W]\subseteq cl(V))
\end{array} $\\
$\begin{array}{c}
\Rightarrow (\exists U\in eR(X,x))(\exists W\in eR(X,y))(f[U]\cap f[W] = \emptyset)
\end{array}$\\
$\begin{array}{c}
\Rightarrow (\exists U\in eR(X,x))(\exists W\in eR(X,y))(U \cap V = \emptyset).
\end{array}$
\end{pf}

\begin{Defn}
A space $X$ is called to be $eR$-compact if every cover of $X$ by $e$-regular sets has a finite subcover.
\end{Defn}

\begin{thm}
Let $f: X \rightarrow Y$ be a function having an $er$-graph $G(f)$, then $f[K]$ is $\theta$-closed in $Y$ for all $eR$-compact relative to $X$ subset $K$.
\end{thm}

\begin{pf}
Let $K$ be $eR$-compact relative to $X$ and $y\notin f[K].$ \newline
$\left.
\begin{array}{r}
y\notin f[K]\Rightarrow (\forall x\in K)(\left( x,y\right) \notin G(f)) \\
G(f) \mbox{ is $er$-graph}
\end{array}%
\right\} \overset{\text{Lemma \ref{10}}}{\Rightarrow} $\newline
$\begin{array}{c}
\Rightarrow (\exists U_{x}\in eR(X,x))(\exists V_{x}\in O(Y,y))(f[U_{x}]\cap cl( V_{x})=\emptyset )
\end{array}$\\
$\left.
\begin{array}{r}
\Rightarrow (\{U_{x}|x\in K\}\subseteq eR(X))(K\subseteq \bigcup \{U_{x}|x\in K\})%
\\
K\mbox{ is }eR\mbox{-compact relative to }X%
\end{array}%
\right\} \Rightarrow $\newline
$\left.
\begin{array}{r}
\Rightarrow (\exists K^{\ast }\subseteq K)(|K^{\ast }|<\aleph_{0})\newline
(K\subseteq \bigcup \{U_{x}|x\in K^{\ast }\}) \\
V:=\underset{x\in K^{\ast }}{\bigcap}V_{x}\in O(Y,y)%
\end{array}%
\right\} \Rightarrow $\newline
$\left.
\begin{array}{c}
\Rightarrow \left( V\in O(Y,y)\right) \left( f[K]\cap V\subseteq \left(
\underset{x\in K^{\ast }}{\bigcup }f\left[U_{x}\right]  \right) 
\cap cl(V) \subseteq 
\underset{x\in K^{\ast }}{\bigcup}\left(f\left[U_{x}\right]   
\cap cl(V) \right)=\emptyset \right) \newline
\end{array}%
\newline
\right. $\newline
$\left.
\begin{array}{c}
\Rightarrow y\notin cl_{\theta }(f[K]).\newline
\end{array} 
\newline
\right. $ 
\end{pf}

\section{Conclusion}
This paper is concerned with the concept of $eR$-continuity and weakly $eR$-continuity defined by utilizing the notion of $e$-regular sets. It turns out that both $eR$-continuous and weakly $eR$-continuous functions are stronger than weakly $e$-continuous functions, as will be seen in Remark \ref{q}. We believe that this paper will pave the way for future studies relevant to continuity and convergence etc. known from functional analysis.

%%%%%%%%%%%%%%%%%%%%%%%%%%%%%%%%%%%%%%%%%%%%%%%%%%%%%%%%%%%%%%%%%%%

\ack The author would like to thank to anonymous reviewer for his/her valuable suggestions and comments, which improved this study.

%%%%%%%%%%%%%%%%%%%%%%%%%%%%%%%%%%%%%%%%%%%%%%%%%%%%%%%%%%%%%%%%%%%
%\section{Bibliography}

\end{document}